\documentclass{article}
\usepackage{comment,fullpage,amssymb,amsmath,amsthm,enumerate,hyperref,mathrsfs,tikz-cd,bbm,bm,adjustbox}
\usepackage[utf8]{inputenc}
\usepackage[style=alphabetic, maxnames=10]{biblatex}
\theoremstyle{plain}
\newtheorem{theorem}{Theorem}[section]

\newtheorem{corollary}[theorem]{Corollary}

\newtheorem{lemma}[theorem]{Lemma}

\theoremstyle{definition}

\theoremstyle{remark}

\newtheorem{remark}[theorem]{Remark}

\title{Large sieve inequalities with power moduli and Waring's problem}
\author{Stephan Baier and Sean B. Lynch}

\newcommand{\Addresses}{{
  \bigskip
  \footnotesize

  \textsc{Stephan Baier, Department of Mathematics, Ramakrishna Mission Vivekananda Educational Research Institute, G. T. Road, PO Belur Math, Howrah, West Bengal 711202, India}\par\nopagebreak
  \textit{E-mail address}: \texttt{stephanbaier2017@gmail.com}

  \medskip

  \textsc{Sean B. Lynch, Institute for Theoretical Sciences, Westlake University in Hangzhou, Zhejiang, China}\par\nopagebreak
  \textit{E-mail address}: \texttt{sblynch@westlake.edu.cn}

}}

\begin{document}

\maketitle

\begin{abstract}
We improve the large sieve inequality with $k\textsuperscript{th}$-power moduli, for all $k\ge 4$. Our method relates these inequalities to a variant of Waring's problem with restricted $k\textsuperscript{th}$-powers. Firstly, we input a classical divisor bound on the number of representations of a positive integer as a sum of two $k\textsuperscript{th}$-powers. Secondly, we apply Marmon's bound on the number of representations of a positive integer as a sum of four $k\textsuperscript{th}$-powers. Thirdly, we use Wooley's Vinogradov mean value theorem with arbitrary weights. Lastly, we state a conditional result, based on the conjectural Hardy-Littlewood formula for the number of representations of a large positive integer as a sum of $k+1$ $k\textsuperscript{th}$-powers.
\end{abstract}

\section{Introduction}

Let $M$ be an integer, let $N$ be a positive integer, and let $\{z_n\}$ be a sequence of complex numbers. Let $\mathcal{S}$ be a finite set of real numbers, and let $\delta\le 1$ be a positive real number. We write $\lVert \alpha\rVert$ for the distance from a real number $\alpha$ to its nearest integer, and we also write $e(\alpha)=e^{2\pi i\alpha}$. Moreover, we say that the set $\mathcal{S}$ is \textit{$\delta$-spaced modulo 1} if $\lVert \alpha-\beta\rVert\ge\delta$ whenever $\alpha$ and $\beta$ are distinct members of $\mathcal{S}$. The general large sieve inequality says that if $\mathcal{S}$ is $\delta$-spaced modulo 1, then
\begin{align}\label{general_lsi}
    \sum_{\alpha\in \mathcal{S}}\left|\sum_{n=M+1}^{M+N} z_n e(n\alpha)\right|^2\le (\delta^{-1}+N-1)\sum_{n=M+1}^{M+N}|z_n|^2.
\end{align}
Montgomery wrote a lovely expository article on this result \cite{Montgomery78}. 

Let $Q$ and $k$ be positive integers. Then we take the set 
\begin{align}\label{set_farey}
    \mathcal{S}_k(Q)=\left\{\text{$a/q^k: a $ and $q$ are coprime positive integers with $a\le q^k\le Q^k$}\right\}
\end{align}
so that
$$\sum_{\alpha\in \mathcal{S}_k(Q)}\left|\sum_{n=M+1}^{M+N} z_n e(n\alpha)\right|^2=\sum_{q=1}^Q\sum_{\substack{a=1\\ \gcd(a,q)=1}}^{q^k}\left|\sum_{n=M+1}^{M+N} z_n e(na/q^k)\right|^2.$$
This paper is about large sieve inequalities with $k\textsuperscript{th}$-power moduli. These are inequalities of the form 
$$\sum_{q=1}^Q\sum_{\substack{a=1\\ \gcd(a,q)=1}}^{q^k}\left|\sum_{n=M+1}^{M+N} z_n e(n a/q^k)\right|^2\le \Delta_k(Q,N) \sum_{n=M+1}^{M+N}|z_n|^2,$$
where $\Delta_k(Q,N)$ is independent of $M$ and $\{z_n\}$. We now let $\Delta_k(Q,N)$ be the infimum of such constants so that our large sieve inequalities are just upper bounds on $\Delta_k(Q,N)$.

 In the first paper on these bounds, Zhao observed two bounds that trivially follow from the general large sieve inequality \cite{LiangyiZhao2004}. Firstly, since the set $\mathcal{S}_k(Q)$ is $Q^{-2k}$-spaced modulo 1, we have
\begin{align}\label{proof_first_trivial}
    \Delta_k(Q,N)\le Q^{2k}+N-1.
\end{align}
Secondly, we can apply the general large sieve inequality to the sum over $a$, and then sum over $q$. This gives
\begin{align}\label{proof_second_trivial}
    \Delta_k(Q,N)\le \sum_{q=1}^Q(q^k+N-1)\le Q(Q^k+N-1).
\end{align}
If $N\gg Q^{2k}$, then bound (\ref{proof_first_trivial}) gives $\Delta_k(Q,N)\ll N$. If $N\ll Q^k$, then bound (\ref{proof_second_trivial}) gives $\Delta_k(Q,N)\ll Q^{k+1}$. Combining these, we have
\begin{align}\label{combined_trivial}
    N\gg Q^{2k} \mbox{ or } N\ll Q^k \implies \Delta_k(Q,N)\ll Q^{k+1}+N.
\end{align}

It's not hard to see that 
\begin{align}\label{lower_bound}
    \Delta_k(Q,N)\gg Q^{k+1}+N.
\end{align}
So, in particular, (\ref{combined_trivial}) is sharp. If $z_{M+1}=1$ and every other $z_n=0$, then we compute that
$$\sum_{\alpha\in \mathcal{S}_k(Q)}\left|\sum_{n=M+1}^{M+N} z_n e(n \alpha)\right|^2=|\mathcal{S}_k(Q)|\sum_{n=M+1}^{M+N}|z_n|^2.$$
Therefore
\begin{align}\label{lower_bound_part_1}
    \Delta_k(Q,N)\ge |\mathcal{S}_k(Q)|=\sum_{q=1}^Q \varphi(q^k)=\sum_{q=1}^Q q^{k-1}\varphi(q)\gg Q^{k+1}.
\end{align}
If each $z_n=e(-n\alpha_0)$, for a fixed $\alpha_0\in\mathcal{S}_k(Q)$, then we have
$$\sum_{\alpha\in \mathcal{S}_k(Q)}\left|\sum_{n=M+1}^{M+N} z_n e(n \alpha)\right|^2\ge \left|\sum_{n=M+1}^{M+N} z_n e(n \alpha_0)\right|^2=N\sum_{n=M+1}^{M+N}|z_n|^2.$$
Therefore
\begin{align}\label{lower_bound_part_2}
    \Delta_k(Q,N)\ge N.
\end{align}
Combining bounds (\ref{lower_bound_part_1}) and (\ref{lower_bound_part_2}), we get bound (\ref{lower_bound}).

In the case $k=1$, bounds (\ref{proof_first_trivial}) and (\ref{lower_bound}) collectively show that $\Delta_1(Q,N)\asymp Q^2+N$. Otherwise, we call
\begin{align}\label{nontrivial_range}
    Q^k<N<Q^{2k}
\end{align}
the \textit{nontrivial range}. This is where many authors \cite{LiangyiZhao2004,BaierZhao05,Baier06,BaierZhao08,Halupczok2012,Halupczok2015,Halupczok2018,Halupczok2020,Munsch2021,BakerMunschShparlinski2022,mcgrath2022asymmetric} beat the trivial bounds and worked toward Zhao's conjecture that   
\begin{align*}
    \Delta_k(Q,N)\ll Q^\varepsilon(Q^{k+1}+N).
\end{align*}
 Baier, Lynch and Zhao showed that the factor of $Q^\varepsilon$ cannot be removed when $k=2$ \cite{baier_lynch_zhao_2019}. Their method extends to all $k\ge 2$.

Our method of attacking the nontrivial range is simple, effective and can easily be applied to more general large sieve inequalities with polynomial moduli. We choose to illustrate our method in this restricted setting of $k\textsuperscript{th}$-power moduli for three reasons. Firstly, these large sieve inequalities have received a lot of attention, and we achieve savings in substantial ranges for all $k\ge 4$. Secondly, they have many (actualised and potential) applications, as listed in the introduction to \cite{BakerMunschShparlinski2022}. Lastly, our method relates these large sieve inequalities to Waring's problem with $k\textsuperscript{th}$-powers. 

For an introduction to Vinogradov's mean value theorem and its application to Waring's problem, look to Vaughan's book \cite[Ch.~5]{Vaughan97}. In the case $k=3$, Wooley proved a long standing conjectured bound on Vinogradov mean values \cite[Thm.~1.1]{Wooley2016}. Bourgain, Demeter and Guth were the first to prove the conjectured bound for all $k\ge 3$ \cite[Thm.~1.1]{Bourgain2016}. We make use of Wooley's generalised version \cite[Thm.~1.1]{Wooley2018} of Bourgain, Demeter and Guth's bound. See \cite{cook2022decoupling} for an interesting comparison between the method of Wooley and that of Bourgain, Demeter and Guth.

\bigskip

\textbf{Acknowledgements:} The first author thanks the Ramakrishna Mission Vivekananda Education and Research Institute for an excellent work environment. The second author thanks the University of New South Wales and the Australian Research Council (DP220102861) for financially supporting this work. We both thank Marc Munsch, Igor Shparlinski and anonymous referees for their comments on the first draft of this paper.

\section{Results}

Let $Q,N,n$ and $k,s,t$ be positive integers. In light of the introduction, we restrict to the nontrivial range
\begin{align*}
    Q^{k}< N < Q^{2k}.
\end{align*}
Let $R_{k,s}(n)$ be the number of $s$-tuples of positive integers $(n_1,\ldots,n_s)$ satisfying $n=n_1^k+\ldots+n_s^k$. The forthcoming implied constants are independent of $Q,N,n$. They may depend on $k,s,t$ or on an arbitrarily small positive real number $\varepsilon$. 

We pass to counting Farey fractions in short intervals in Section \ref{section_ls_to_counting_fractions}. This is standard practice. We then lay the foundations of our new method in Section \ref{section_counting_fractions_to_waring}, where we also prove the following two lemmata.

\begin{lemma}\label{waring_to_lsi}
Suppose that, uniformly in positive integers $n\le sQ^k$, we have the bound
\begin{align*}
    R_{k,s}(n) \ll C_{\varepsilon,k,s}(Q).
\end{align*}
Then we have 
    \begin{align*}
        \Delta_k(Q,N)\ll NQ^\varepsilon \left(C_{\varepsilon,k,s}(Q)\cdot N^{-1}Q^{2k}  \right)^{1/s}.
    \end{align*}
\end{lemma}

\begin{lemma}\label{mvt_to_lsi}
    Suppose that, uniformly in subsets $\mathcal{A}$ of positive integers $\le Q$, we have the bound
    \begin{align}\label{mvt_estimate}
        \int_0^1 \left|\sum_{n\in\mathcal{A}} e(\alpha n^k)\right|^{2s} d\alpha\ll D_{\varepsilon,k,s}(Q)\cdot |\mathcal{A}|^s.
    \end{align}
Then we have 
\begin{align*}
    \Delta_k(Q,N)\ll NQ^\varepsilon\left(D_{\varepsilon,k,s}(Q)\cdot N^{-1}Q^{2k}  \right)^{1/s}.
\end{align*}
\end{lemma}

In Section \ref{section_sup_proof}, we recall an elementary divisor bound for $R_{k,2}(n)$. Running this through Lemma \ref{waring_to_lsi} with $s=2$ recovers Baier and Zhao's large sieve inequality \cite[Thm.~1.1]{BaierZhao05}. We state this large sieve inequality in part (i) of our first theorem. 

Also in Section \ref{section_sup_proof}, we recall Marmon's non-elementary bound on $R_{k,4}(n)$ \cite[Thm.~1.3]{Marmon2010}. Marc Munsch suggested that we insert this into Lemma \ref{waring_to_lsi} with $s=4$ so as to recover Baker, Munsch and Shparlinski's large sieve inequality \cite[Thm.~1.4]{BakerMunschShparlinski2022}. We state this large sieve inequality in part (ii) of our first theorem. In fact, we state the inequality for all $k\ge 3$, while Baker, Munsch and Shparlinski left a void when $k=4$. McGrath \cite[Cor.~1.4]{mcgrath2022asymmetric} filled this void, albeit with a weaker bound than the general formula in $k$ would suggest. Our result corrects this, filling the void with the appropriate bound.  

\begin{theorem}\label{thm_sup_s=2}
\begin{enumerate}[(i)]
    \item For all $k\ge 2$, we have the large sieve inequality
\begin{align*}
    \Delta_k(Q,N)\ll N^{\frac1{2}}Q^{k+\varepsilon}.
\end{align*}
\item For all $k\ge 3$, we have the large sieve inequality
\begin{align*}
    \Delta_k(Q,N)\ll N^{\frac{3}{4}}Q^{\frac{k}{2}+\frac{1}{4}+\frac{1}{2\sqrt{k}}+\varepsilon}.
\end{align*}
\end{enumerate}
\end{theorem}

Finally, in Section \ref{section_moment_proof}, we invoke a recent and high-powered theorem of Wooley \cite[Thm.~1.1]{Wooley2018}. We end up with a family, indexed by positive integers $t\le k$, of mean value estimates of the form (\ref{mvt_estimate}). Running this through Lemma \ref{mvt_to_lsi} with $s=\frac1{2}t(t+1)$ results in our second theorem.

\begin{theorem}\label{thm_second_moment}
    For any positive integer $t\le k$, we have the large sieve inequality
\begin{align}\label{lsi_in_t}
    \Delta_k(Q,N)\ll N^{1-\frac{2}{t(t+1)}}Q^{1+\frac{4k-2t}{t(t+1)}+\varepsilon}.
\end{align}
\end{theorem}

\begin{remark}
Write $N=Q^{\lambda}$ for a real $\lambda$ with $k<\lambda<2k$. Then the right-hand side of (\ref{lsi_in_t}) becomes
\begin{align*}
    Q^{f_{k,t}(\lambda)+\varepsilon},
\end{align*}
where
\begin{align*}
    f_{k,t}(\lambda)=\lambda\left(1-\frac{2}{t(t+1)}\right)+1+\frac{4k-2t}{t(t+1)}=\lambda+1+\frac{2(2k-\lambda)}{t(t+1)}-\frac{2}{t+1}.
\end{align*}
We minimise $f_{k,t}(\lambda)$ by taking 
$$t=\min\{\lceil2(2k-\lambda)\rceil,k\}.$$
Incidentally, this choice of $t$ makes $f_{k,t}(\lambda)$ continuous in $\lambda$.
\end{remark}

As we'll see in the next section, Theorem \ref{thm_second_moment} with $t=k$ is particularly powerful.

\begin{corollary}\label{cor_main_bound}
    We have the large sieve inequality
\begin{align*}
    \Delta_k(Q,N)\ll N^{1-\frac{2}{k(k+1)}}Q^{1+\frac{2}{k+1}+\varepsilon}.
\end{align*}
\end{corollary}

To discuss a conditional result, we assume that $k\ge 3$ and $s\ge k+1$ throughout the rest of this section. Recall the Hardy-Littlewood singular series
\begin{align*}
    \mathfrak{S}_{k,s}(n)=\sum_{q=1}^\infty\sum_{\substack{a=1\\ \gcd(a,q)=1}}^q \left(\frac1{q} \sum_{r=1}^q e(a r^k/q)\right)^s e(-na/q).
\end{align*}
Then Hardy and Littlewood proved that whenever $s$ is sufficiently large relative to $k$, we have the formula
\begin{align}\label{hardy-littlewood_asymptotic}
    R_{k,s}(n)=\frac{\Gamma(1+\frac1{k})^s}{\Gamma(\frac{s}{k})}\mathfrak{S}_{k,s}(n)n^{\frac{s}{k}-1}+o(n^{\frac{s}{k}-1}).
\end{align}
See \cite[Ch.~2]{Vaughan97} for a modern proof of this formula, with $s>2^k$ and a power saving in the error term. We are often interested in obtaining $\mathfrak{S}_{k,s}(n)\gg 1$ so that the first term on the right-hand side is not swallowed by the error, thereby making (\ref{hardy-littlewood_asymptotic}) a fully-fledged asymptotic formula. This lower-bound on the singular series can be subtle, as covered in \cite[\S.~4.5]{Vaughan97}. Fortunately, for our application, the only bound we need on the singular series is a modest upper-bound. Indeed, we have \cite[Thm.~4.3, p.~49]{Vaughan97}
$$\mathfrak{S}_{k,s}(n)\ll n^\varepsilon.$$
Therefore, if formula (\ref{hardy-littlewood_asymptotic}) holds, then we have 
\begin{align}\label{typical_waring_bound}
    R_{k,s}(n)\ll n^{\frac{s}{k}-1+\varepsilon}.
\end{align}
If bound (\ref{typical_waring_bound}) holds, then it follows from Lemma \ref{waring_to_lsi} that
\begin{align}\label{conditional_lsi}
    \Delta_k(Q,N)\ll N^{1-\frac{1}{s}}Q^{1+\frac{k}{s}+\varepsilon}.
\end{align}
Rewriting this large sieve bound as $NQ^{1+\varepsilon}(N^{-1}Q^k)^{1/s}$ and recalling that $N> Q^k$, it's optimal to choose $s$ as small as possible. As luck would have it, this runs with -- not against -- the difficulty in Waring's problem! The Hardy-Littlewood asymptotic formula (\ref{hardy-littlewood_asymptotic}) is conjectured to hold with $s=k+1$ whenever $k\ge 3$. If this famous conjecture holds, then large sieve bound (\ref{conditional_lsi}) holds with $s=k+1$ whenever $k\ge 3$. 

\section{Comparison with previous results}

Let $Q,N,k$ be positive integers. The forthcoming implied constants are independent of $Q$ and $N$, but they may depend on $k$ or on an arbitrarily small positive real number $\varepsilon$. 

The standard way of comparing bounds is to write $N=Q^{\lambda}$ for a real parameter $\lambda$. Restricting to the nontrivial range (\ref{nontrivial_range}), we assume that $k<\lambda<2k$. It is also customary to introduce clean notation for some functions of $k$ that appear in exponents. Put
\begin{align*}
    \kappa=\frac1{2^{k-1}} \qquad \mbox{and} \qquad \omega=\frac1{(k-1)(k-2)+2}.
\end{align*}

\subsection{\texorpdfstring{$k\ge 5$}{}}
In this subsection, we assume that $k\ge 5$. Here is the list of bounds with which we have to contend:\\
First trivial (\ref{proof_first_trivial})
\begin{align}\label{bound_first_trivial}
    \Delta_k(Q,N)\ll Q^{2k}+N,
\end{align}
Second trivial (\ref{proof_second_trivial})
\begin{align}\label{bound_second_trivial}
    \Delta_k(Q,N)\ll Q(Q^k+N),
\end{align}
Zhao \cite{LiangyiZhao2004}
\begin{align}\label{bound_zhao}
    \Delta_k(Q,N)\ll (Q^{k+1}+NQ^{1-\kappa}+N^{1-\kappa}Q^{1+k\kappa})Q^\varepsilon,
\end{align}
Baier and Zhao \cite{BaierZhao05}
\begin{align}\label{bound_baier-zhao}
    \Delta_k(Q,N)\ll (Q^{k+1}+N+N^{\frac1{2}}Q^{k})Q^{\varepsilon},
\end{align}
Halupczok \cite{Halupczok2012,Halupczok2015,Halupczok2018,Halupczok2020}
\begin{align}\label{bound_halupczok}
    \Delta_k(Q,N)\ll \left(Q^{k+1}+\min\left\{N^{1-\omega}Q^{1+(2k-1)\omega},NQ^{1-\frac1{k(k-1)}}+N^{1-\frac1{k(k-1)}}Q^{\frac{k}{k-1}}\right\}\right)Q^{\varepsilon},
\end{align}
Munsch \cite{Munsch2021}
\begin{align}\label{bound_munsch}
    \Delta_k(Q,N)\ll N^{1-\frac1{k(k+1)}}Q^{1+\frac{1}{k+1}+\varepsilon},
\end{align}
Baker, Munsch and Shparlinski \cite{BakerMunschShparlinski2022}
\begin{align}\label{bound_baker-munsch-shparlinski}
    \Delta_k(Q,N)\ll NQ^{\frac1{2}}+N^{\frac{3}{4}}Q^{\frac{k}{2}+\frac1{4}+\frac{1}{2\sqrt{k}}+\varepsilon}.
\end{align}

Corollary \ref{cor_main_bound} is our main bound. It is sharper than all of the previous bounds for $\lambda<\lambda_0$, where $\lambda_0=2k-3+O(k^{-\frac1{2}})$ and $2k-3<\lambda_0<2k-1$. As $k$ tends to infinity, our main bound tends to improve $100\%$ of the nontrivial range. In particular, when $k=5$, it already improves more than $69\%$ of the nontrivial range. This percentage is more than $86\%$ when $k=16$, and it is more than $97\%$ when $k=100$.

For the remaining range $\lambda\ge\lambda_0$, (\ref{bound_baier-zhao}) and (\ref{bound_baker-munsch-shparlinski}) are the strongest of the previous bounds. Theorem \ref{thm_sup_s=2} recovers both of these, and the proofs we give are simpler than the previous proofs. Most notably, Marc Munsch observed that we may completely bypass Baker, Munsch and Shparlinski's route through asymmetric additive energy; they use a consequence of Marmon's three-variable bound \cite[Thm.~1.4]{Marmon2010} to prove (\ref{bound_baker-munsch-shparlinski}). 

By applying Theorem \ref{thm_second_moment} with smaller values of $t$, we increase the range in which we improve on these two previous bounds. Indeed, Theorem \ref{thm_second_moment} with $t=3$ (or even $t=4$) beats (\ref{bound_baier-zhao}) for $\lambda<2k-\frac{3}{2}$. The choice of $t$ that maximises our range of improvement on (\ref{bound_baker-munsch-shparlinski}) changes with $k$. When $5\le k<16$, we take $t=3$ to improve on (\ref{bound_baker-munsch-shparlinski}) in the range $\lambda<2k-3+\frac{6}{\sqrt{k}}$. When $16\le k<100$, we take $t=4$ to improve on (\ref{bound_baker-munsch-shparlinski}) in the range $\lambda<2k-\frac{7}{3}+\frac{10}{3\sqrt{k}}$. Finally, when $k\ge 100$, we take $t=5$ to improve on (\ref{bound_baker-munsch-shparlinski}) in the range $\lambda<2k-\frac{25}{11}+\frac{30}{11\sqrt{k}}$.

Now let
$$\lambda_1=\begin{cases}
    2k-\frac{3}{2} & \mbox{ if } 5\le k<16\\
    2k-\frac{7}{3}+\frac{10}{3\sqrt{k}} & \mbox{ if } 16\le k<100\\
    2k-\frac{25}{11}+\frac{30}{11\sqrt{k}} & \mbox{ if } k\ge 100.
\end{cases}$$
Then Theorem \ref{thm_second_moment}, taken as one whole, beats all of the previous bounds for $\lambda<\lambda_1$. So Theorem \ref{thm_second_moment} improves $70\%$ of the nontrivial range when $k=5$. This percentage is more than $90\%$ when $k=16$, and it is exactly $98\%$ when $k=100$.

\subsection{\texorpdfstring{$k=4$}{}}

We have to contest these previous bounds: (\ref{bound_first_trivial}), (\ref{bound_second_trivial}), (\ref{bound_zhao}), (\ref{bound_baier-zhao}), (\ref{bound_halupczok}), (\ref{bound_munsch}). Indeed, (\ref{bound_baker-munsch-shparlinski}) is missing. 

Corollary \ref{cor_main_bound} beats all of these previous bounds when $\frac{21}{4}<\lambda<\frac{13}{2}$. Hence it improves on the previous bounds in over $31\%$ of the nontrivial range. Smaller choices of $t$ in Theorem \ref{thm_second_moment} do not increase our range of improvement. While Theorem \ref{thm_sup_s=2} part (ii) does bring forth (\ref{bound_baker-munsch-shparlinski}), it also fails to increase our range of improvement.  

In the range $\lambda<\frac{21}{4}$, (\ref{bound_zhao}) is the strongest; we cannot even match this. In the range $\lambda\ge \frac{13}{2}$, (\ref{bound_baier-zhao}) is the strongest; Theorem \ref{thm_sup_s=2} part (i) recovers this.

\subsection{\texorpdfstring{$k=3$}{}}

In the case of cubes, Baier and Zhao proved the large sieve inequality \cite{BaierZhao05}
\begin{align}\label{bound_baier-zhao_cubic}
    \Delta_3(Q,N)\ll (Q^4+N^{\frac{9}{10}}Q^{\frac{6}{5}}+NQ^{\frac{6}{7}})Q^\varepsilon.
\end{align}
This is the strongest bound in the range $\lambda\le \frac{25}{7}$. Then we see that (\ref{bound_zhao}) is the strongest in the range $\frac{25}{7}\le \lambda\le \frac{9}{2}$, and (\ref{bound_baier-zhao}) is the strongest in the range $\lambda\ge \frac{9}{2}$. While we offer no improvement, Theorem \ref{thm_sup_s=2} part (i) does match (\ref{bound_baier-zhao}).

\subsection{\texorpdfstring{$k=2$}{}}

In the case of squares, Baier and Zhao gave the large sieve inequality \cite{BaierZhao08}
$$\Delta_2(Q,N)\ll (Q^3+N+\min\{NQ^{\frac1{2}},N^{\frac1{2}}Q^2\})Q^\varepsilon.$$
This is the state of the art. While we offer no improvement, Theorem \ref{thm_sup_s=2} part (i) does match it when $\lambda\ge 3$.

\section{From the large sieve to Farey fractions in short intervals}\label{section_ls_to_counting_fractions}

Let $Q,N$ and $k,s,t$ be positive integers. We also let $M$ be an integer, and let $\{z_n\}$ be a sequence of complex numbers. The forthcoming implied constants are independent of $Q,M,N$ and $\{z_n\}$. But they may depend on $k,s,t$ or on an arbitrarily small positive real number $\varepsilon$. We assume that
\begin{align*}
    Q^k\le N\ll Q^{2k}.
\end{align*}

Dividing the summation over $q$ into dyadic intervals gives
\begin{align*}
    \sum_{q=1}^Q\sum_{\substack{a=1\\ \gcd(a,q)=1}}^{q^k}\left|\sum_{n=M+1}^{M+N} z_n e(n a/q^k)\right|^2\ll Q^\varepsilon \sup_{1\le x\le Q} \sum_{\frac{x}{2}<q\le x}\sum_{\substack{a=1\\ \gcd(a,q)=1}}^{q^k}\left|\sum_{n=M+1}^{M+N} z_n e(n a/q^k)\right|^2.
\end{align*}
Let $x$ be a real number in the interval $1\le x\le Q$ and write 
$$\mathcal{A}_k(x):=\left\{\frac{a}{q^k}: 1\le a\le q^k,\ \frac{x}{2}<q\le x \ \mbox{and}\, \gcd(a,q)=1\right\}$$
so that
$$\sum_{\frac{x}{2}<q\le x}\sum_{\substack{a=1\\ \gcd(a,q)=1}}^{q^k}\left|\sum_{n=M+1}^{M+N} z_n e( na/q^k)\right|^2=\sum_{\alpha\in\mathcal{A}_k(x)}\left|\sum_{n=M+1}^{M+N} z_n e(n\alpha)\right|^2.$$
 We want to hit this with the general large sieve inequality (\ref{general_lsi}), leading us to consider the modulo 1 spacing of the set $\mathcal{A}_k(x)$. Let $\delta=\delta(x)\le x^{-k}$ be a positive real number. Note that $\min\mathcal{A}_k(x)\ge x^{-k}$ and that $\max\mathcal{A}_k(x)\le 1$. Then, for every $\alpha,\beta\in\mathcal{A}_k(x)$, we have $\lVert\alpha-\beta\rVert<\delta$ if and only if $|\alpha-\beta|<\delta$. Hence, for every $\alpha_0\in\mathcal{A}_k(x)$, we have
$$\mathcal{A}_k(x,\delta,\alpha_0):=\{\alpha\in\mathcal{A}_k(x):\lVert\alpha-\alpha_0\rVert<\delta\}=\{\alpha\in\mathcal{A}_k(x):|\alpha-\alpha_0|<\delta\}.$$
Now observe that there exists a partition of the set $\mathcal{A}_k(x)$ into 
$$\sup_{\alpha_0\in\mathcal{A}_k(x)}|\mathcal{A}_k(x,\delta,\alpha_0)|$$
parts so that each part is $\delta$-spaced modulo 1. Thus, for any one of these parts $\mathcal{P}$, say, the general large sieve inequality (\ref{general_lsi}) gives
$$\sum_{\alpha\in\mathcal{P}}\left|\sum_{n=M+1}^{M+N} z_n e(n \alpha)\right|^2\le (N-1+\delta^{-1})\sum_{n=M+1}^{M+N} |z_n|^2.$$
This then yields the bound
$$\sum_{\alpha\in\mathcal{A}_k(x)}\left|\sum_{n=M+1}^{M+N} z_n e( n\alpha)\right|^2\le\left(\sup_{\alpha_0\in\mathcal{A}_k(x)}|\mathcal{A}_k(x,\delta,\alpha_0)|\right)(N-1+\delta^{-1})\sum_{n=M+1}^{M+N} |z_n|^2.$$
Putting everything together, we may take
\begin{align*}
    \Delta_k(Q,N)\ll Q^\varepsilon\sup_{1\le x\le Q}\ \inf_{0<\delta\le x^{-k}}\left((N+\delta^{-1})\sup_{\alpha_0\in\mathcal{A}_k(x)}|\mathcal{A}_k(x,\delta,\alpha_0)|\right).
\end{align*}
Recall that $x^k\le Q^k\le N$. Therefore $\frac1{2N}\le \frac{1}{2x^k}$, and so we may take 
$$\delta=\frac1{2N}$$ 
in the infimum. Hence, with this choice of $\delta$, we have 
\begin{align}\label{standard_opening}
    \Delta_k(Q,N)\ll Q^\varepsilon N\sup_{1\le x\le Q}\sup_{\alpha_0\in\mathcal{A}_k(x)}|\mathcal{A}_k(x,\delta,\alpha_0)|.
\end{align}

\section{From Farey fractions in short intervals to Waring's problem}\label{section_counting_fractions_to_waring}

We keep the notation of the previous section. Recall that $x$ is a real number in the interval $1\le x\le Q$ and that we have $$\delta=\frac1{2N}\le\frac{1}{2x^k}.$$
Let $\alpha_0\in \mathcal{A}_k(x)$. Then we write $a_0$ and $q_0$ for the positive integers satisfying
$$\alpha_0=\frac{a_0}{q_0^k},\qquad a_0\le q_0^k, \qquad \frac{x}{2}<q_0\le x \qquad \mbox{and} \qquad \gcd(a_0,q_0)=1.$$
Let $\mathcal{A}^*_k(x,\delta,a_0,q_0)$ be the set of positive integers $q$ in the dyadic interval $\frac{x}{2}<q\le x$ for which there exists an integer $b$ satisfying $|b|<\delta x^{2k}$ and $a_0q^k\equiv b\bmod q_0^k$.

 Given any $\alpha\in\mathcal{A}_k(x,\delta,\alpha_0)$, we write $\alpha=a/q^k$ in reduced form like we did for $\alpha_0$. Then we rewrite the condition $|\alpha-\alpha_0|<\delta$ as $|a_0q^k-aq_0^k|<\delta q^k q_0^k\le \delta x^{2k}$. Therefore $\alpha\mapsto q$ defines a function $$\mathcal{A}_k(x,\delta,\alpha_0)\hookrightarrow \mathcal{A}^*_k(x,\delta,a_0,q_0).$$ 
 Since $\delta$ is sufficiently small, this function is injective. To see this, take two fractions $\alpha,\beta\in\mathcal{A}_k(x,\delta,\alpha_0)$ with the same reduced denominator $q^k$. Now observe that $|\alpha-\beta|\le |\alpha-\alpha_0|+|\alpha_0-\beta|<2\delta$. Then we see that $\alpha=\beta$, because otherwise we would get the contradictory bound $|\alpha-\beta|\ge q^{-k}\ge x^{-k}\ge 2\delta$. Hence
\begin{align}\label{injection_bound}
    |\mathcal{A}_k(x,\delta,\alpha_0)|\le |\mathcal{A}^*_k(x,\delta,a_0,q_0)|.
\end{align}

Now we let $\mathcal{B}^*_{k,s}(x,\delta,a_0,q_0)$ be the set of positive integers $q$ in the interval $s(\frac{x}{2})^k<q\le sx^k$ for which there exists an integer $b$ satisfying
\begin{align}\label{b_residue_condition}
    |b|<s\delta x^{2k} \qquad \mbox{and} \qquad a_0q\equiv b\bmod q_0^k.
\end{align}
Evidently, if $q_1,\ldots,q_s\in\mathcal{A}^*_k(x,\delta,a_0,q_0)$, then $q_1^k+\ldots+q_s^k\in\mathcal{B}^*_{k,s}(x,\delta,a_0,q_0)$. This leads us to consider, for each $q\in \mathcal{B}^*_{k,s}(x,\delta,a_0,q_0)$, the number $R^*_{k,s}(x,\delta,a_0,q_0;q)$ of $s$-tuples $(q_1,\ldots,q_s)\in \mathcal{A}^*_k(x,\delta,a_0,q_0)^s$ satisfying $q=q_1^k+\ldots+q_s^k$. We have the formula
\begin{align*}
    |\mathcal{A}^*_k(x,\delta,a_0,q_0)|^s= \sum_{q\in\mathcal{B}^*_{k,s}(x,\delta,a_0,q_0)}R^*_{k,s}(x,\delta,a_0,q_0;q).
\end{align*}

We proceed in two ways. First, we take the supremum to get
\begin{align}\label{first_way}
    |\mathcal{A}^*_k(x,\delta,a_0,q_0)|^s\le |\mathcal{B}^*_{k,s}(x,\delta,a_0,q_0)|\sup_{q\in\mathcal{B}^*_{k,s}(x,\delta,a_0,q_0)}R^*_{k,s}(x,\delta,a_0,q_0;q).
\end{align}
Second, we use Cauchy-Schwarz to get
\begin{align}\label{second_way}
    |\mathcal{A}^*_k(x,\delta,a_0,q_0)|^{2s}\le |\mathcal{B}^*_{k,s}(x,\delta,a_0,q_0)| \sum_{q\in\mathcal{B}^*_{k,s}(x,\delta,a_0,q_0)}R^*_{k,s}(x,\delta,a_0,q_0;q)^2.
\end{align}
Either way, we need a bound on $|\mathcal{B}^*_{k,s}(x,\delta,a_0,q_0)|$. 

To count the number of $q\in \mathcal{B}^*_{k,s}(x,\delta,a_0,q_0)$, we first recall condition (\ref{b_residue_condition}), from which we see that there are $O(1+x^{2k}\delta)$ choices for $a_0q\bmod q_0^k$. Then, because $\gcd(a_0,q_0^k)=1$, there are $O(1+x^{2k}\delta)$ choices for $q\bmod q_0^k$. Note that $q\asymp x^k\asymp q_0^k$. Thus, for a fixed residue $r$, there are $O(1)$ choices for $q$ with $q\equiv r\bmod q_0^k$. Putting these estimates together, we have
\begin{align}\label{bound_b_set}
    |\mathcal{B}^*_{k,s}(x,\delta,a_0,q_0)|\ll 1+x^{2k}\delta\ll N^{-1}Q^{2k}.
\end{align}

\subsection{Proof of Lemma \ref{waring_to_lsi}}

Here we continue with the first way. 

For each $q\in \mathcal{B}^*_{k,s}(x,\delta,a_0,q_0)$, we have $R_{k,s}^*(x,\delta,a_0,q_0;q)\le R_{k,s}(q)$. Furthermore, if $q\in \mathcal{B}^*_{k,s}(x,\delta,a_0,q_0)$, then $q\le sx^k\le sQ^k$. Therefore
\begin{align*}
    \sup_{q\in\mathcal{B}^*_{k,s}(x,\delta,a_0,q_0)}R^*_{k,s}(x,\delta,a_0,q_0;q) \le \sup_{q\in\mathcal{B}^*_{k,s}(x,\delta,a_0,q_0)} R_{k,s}(q)\le \sup_{q\le sQ^k} R_{k,s}(q).
\end{align*}
Combining this with bounds (\ref{injection_bound}), (\ref{first_way}) and (\ref{bound_b_set}), we have
\begin{align*}
    |\mathcal{A}_k(x,\delta,\alpha_0)|^s\le |\mathcal{A}^*_k(x,\delta,a_0,q_0)|^s\ll N^{-1}Q^{2k}\sup_{q\le sQ^k} R_{k,s}(q).
\end{align*}
We arrive at Lemma \ref{waring_to_lsi} by running this estimate through (\ref{standard_opening}).   

\subsection{Proof of Lemma \ref{mvt_to_lsi}}

Here we continue with the second way.

Note that 
\begin{align*}
    \sum_{q\in\mathcal{B}^*_{k,s}(x,\delta,a_0,q_0)}R^*_{k,s}(x,\delta,a_0,q_0;q)^2
\end{align*}
counts the number of solutions to the equation
\begin{align}\label{one_equation_system}
    p_1^k+\ldots+p_s^k=q_1^k+\ldots+q_s^k
\end{align}
with each $p_i,q_i\in\mathcal{A}^*_k(x,\delta,a_0,q_0)$. By orthogonality of exponentials, the integral
\begin{align*}
    \int_0^1 \left|\sum_{n\in \mathcal{A}^*_k(x,\delta,a_0,q_0)} e(\alpha n^k)\right|^{2s}d\alpha
\end{align*}
also counts the number of solutions to (\ref{one_equation_system}). 

Suppose that, uniformly in subsets $\mathcal{A}$ of positive integers $\le Q$, we have the bound
\begin{align*}
    \int_0^1 \left|\sum_{n\in \mathcal{A}} e(\alpha n^k)\right|^{2s}d\alpha \ll D_{\varepsilon,k,s}(Q)\cdot|\mathcal{A}|^s.
\end{align*}
Recall that this is the supposition of Lemma \ref{mvt_to_lsi}. If $q\in \mathcal{A}^*_k(x,\delta,a_0,q_0)$, then $q\le x\le Q$. So this mean value estimate holds with $\mathcal{A}=\mathcal{A}^*_k(x,\delta,a_0,q_0)$. Combining this with bounds (\ref{second_way}) and (\ref{bound_b_set}), we have
\begin{align*}
    |\mathcal{A}^*_k(x,\delta,a_0,q_0)|^{2s}\ll N^{-1}Q^{2k}\cdot D_{\varepsilon,k,s}(Q)\cdot|\mathcal{A}^*_k(x,\delta,a_0,q_0)|^s.
\end{align*}
Hence, by bound (\ref{injection_bound}), we have
\begin{align*}
    |\mathcal{A}_k(x,\delta,\alpha_0)|^{s}\le |\mathcal{A}^*_k(x,\delta,a_0,q_0)|^{s}\ll N^{-1}Q^{2k} \cdot D_{\varepsilon,k,s}(Q).
\end{align*}
We arrive at Lemma \ref{mvt_to_lsi} by running this through (\ref{standard_opening}).

\section{Proof of Theorem \ref{thm_sup_s=2}}\label{section_sup_proof}

Let $k$ and $n$ be positive integers. The following implied constants are independent of $n$, but they may depend on $k$ or on an arbitrarily small positive real number $\varepsilon$.

\subsection{Part (i)}

In this subsection, we assume that $k\ge 2$.

From the epsilonic power estimate on the divisor function, we have the bound 
\begin{align}\label{elementary_divisor_bound}
    R_{k,2}(n)\ll n^\varepsilon.
\end{align}
\begin{proof}
Assuming that $k$ is even, we have $R_{k,2}(n)\le R_{2,2}(n)$. So, to finish off this case, we merely apply the classical bound $R_{2,2}(n)\ll n^\varepsilon$ \cite[Ch.~2, \S.~4]{Grosswald85}.

Assume that $k$ is odd, and suppose that $a$ and $b$ are positive integers satisfying $a^k+b^k=n$. Recalling the factorisation $a^k+b^k=(a+b)(a^{k-1}-a^{k-2}b+\cdots -ab^{k-2}+b^{k-1})$, we see that $d=a+b$ divides $n$. Because $k>1$, we have $d\le n<d^k$, and thus also $0<d^k-n\le n^k-n$. Computing $n=a^k+b^k\equiv d^k\bmod a$, we see that $a$ divides $d^k-n$. Hence, as desired, we have
\begin{equation*}
    R_{k,2}(n)\le \sum_{\substack{d\mid n\\
    d\le n<d^k}} \sum_{a\mid (d^k-n)} 1 \ll \sum_{\substack{d\mid n\\
    d\le n<d^k}} n^{\varepsilon/2} \ll n^{\varepsilon}. 
    \qedhere
\end{equation*}
\end{proof}

To get Theorem \ref{thm_sup_s=2} part (i), we merely insert estimate (\ref{elementary_divisor_bound}) into Lemma \ref{waring_to_lsi} with $s=2$.

\subsection{Part (ii)}

In this subsection, we assume that $k\ge 3$.

Using arithmetic algebraic geometry, Marmon derived the bound \cite[Thm.~1.3]{Marmon2010}
\begin{align}\label{marmon_bound}
    R_{k,4}(n)\ll n^{\frac1{k}+\frac{2}{k\sqrt{k}}+\varepsilon}.
\end{align}
The proof is non-elementary, and we don't copy it here. As suggested by Marc Munsch, we get Theorem \ref{thm_sup_s=2} part (ii) by inserting estimate (\ref{marmon_bound}) into Lemma \ref{waring_to_lsi} with $s=4$.

\section{Proof of Theorem \ref{thm_second_moment}}\label{section_moment_proof}

Let $Q$ and $k,s,t$ be positive integers. Let $\mathcal{A}$ be a subset of $\{\text{positive integers $\le Q$}\}$. The forthcoming implied constants are independent of $Q$ and $\mathcal{A}$, but they may depend on $k,s,t$ or on an arbitrarily small positive real number $\varepsilon$.

From Wooley's recent mean value estimate \cite[Thm.~1.1]{Wooley2018}, we get the following family of mean value estimates. For each $t\le k$, we have
\begin{align}\label{mvt}
    \int_0^1 \left|\sum_{n\in\mathcal{A}}e(\alpha n^k)\right|^{t(t+1)} d\alpha \ll Q^{\frac{1}{2}t(t-1)+\varepsilon}|\mathcal{A}|^{\frac{1}{2}t(t+1)}.
\end{align}
\begin{proof}
    We assume that $t\le k$, and we put $s=\frac1{2}t(t+1)$. By orthogonality of exponentials, the integral on the left-hand side of (\ref{mvt}) counts the number of solutions to the equation
    \begin{align}\label{another_one_equation_system}
        m_1^k+\ldots+m_s^k=n_1^k+\ldots+n_s^k
    \end{align}
    with each $m_i,n_i\in\mathcal{A}$. We now inflate this equation to a partial Vinogradov system, a standard technique in the context of Waring's problem. Namely, we consider the system of $t$ equations 
    \begin{equation}\label{partial_vinogradov_system}
  \begin{aligned}
        m_1^k+\ldots+m_s^k&=n_1^k+\ldots+n_s^k\\
        m_1^{t-1}+\ldots+m_s^{t-1}&=n_1^{t-1}+\ldots+n_s^{t-1}\\
         &\vdots\\
        m_1+\ldots+m_s&=n_1+\ldots+n_s
    \end{aligned}
\end{equation}
    with each $m_i,n_i\in\mathcal{A}$. As we explain in Subsection \ref{subsection_explaining_inflation}, when we pass from counting solutions of equation (\ref{another_one_equation_system}) to counting solutions of this partial Vinogradov system (\ref{partial_vinogradov_system}), we incur an inflationary factor of 
    \begin{align*}
    \ll \prod_{j=1}^{t-1} Q^{j}=Q^{\frac1{2}t(t-1)}.
    \end{align*}
    The number of solutions to our partial Vinogradov system is equal to
    \begin{align*}
        \int_{(0,1]^t} \left|\sum_{n\in\mathcal{A}}e(\alpha_1 n+\cdots+\alpha_{t-1}n^{t-1}+\alpha_t n^k)\right|^{2s} d\alpha.
    \end{align*}
    So it suffices to show that this last integral is $\ll Q^\varepsilon|\mathcal{A}|^s$, for which we apply \cite[Thm.~1.1]{Wooley2018}.
\end{proof}

To get Theorem \ref{thm_second_moment}, we merely insert estimate (\ref{mvt}) into Lemma \ref{mvt_to_lsi} with $s=\frac1{2}t(t+1)$.

\subsection{Explaining the inflationary factor}\label{subsection_explaining_inflation}

Here we introduce some quantities that depend on $\mathcal{A}$. Crucially, the implied constants are independent of $\mathcal{A}$. 
For each positive integer $\ell$, let $\mathcal{R}_{k,s}(\mathcal{A},\ell)$ be the set of $s$-tuples $(n_1,\ldots,n_s)\in\mathcal{A}^s$ satisfying $n_1^k+\ldots+n_s^k=\ell$, and let $R_{k,s}(\mathcal{A},\ell)=|\mathcal{R}_{k,s}(\mathcal{A},\ell)|$ be the number of these. Let $$\mathcal{B}_{k,s,t}(\mathcal{A},\ell)=\{(n_1+\ldots+n_s,\ldots,n_1^{t-1}+\ldots+n_s^{t-1}): (n_1,\ldots,n_s)\in\mathcal{R}_{k,s}(\mathcal{A},\ell)\}$$
so that $(n_1,\ldots,n_s)\mapsto (n_1+\ldots+n_s,\ldots,n_1^{t-1}+\ldots+n_s^{t-1})$ defines a function
$$\mathcal{R}_{k,s}(\mathcal{A},\ell)\to \mathcal{B}_{k,s,t}(\mathcal{A},\ell).$$
As per usual, we need to consider the size of the fibres of this function. For each $(\ell_1,\ldots,\ell_{t-1})\in\mathcal{B}_{k,s,t}(\mathcal{A},\ell)$, let $V_{k,s,t}(\mathcal{A},\ell;\ell_1,\ldots,\ell_{t-1})$ be the size of the fibre. Then we have the formula
$$R_{k,s}(\mathcal{A},\ell)=\sum_{(\ell_1,\ldots,\ell_{t-1})\in\mathcal{B}_{k,s,t}(\mathcal{A},\ell)}V_{k,s,t}(\mathcal{A},\ell;\ell_1,\ldots,\ell_{t-1}).$$
By Cauchy-Schwarz, we have
$$R_{k,s}(\mathcal{A},\ell)^2\le |\mathcal{B}_{k,s,t}(\mathcal{A},\ell)|\sum_{(\ell_1,\ldots,\ell_{t-1})\in\mathcal{B}_{k,s,t}(\mathcal{A},\ell)}V_{k,s,t}(\mathcal{A},\ell;\ell_1,\ldots,\ell_{t-1})^2.$$
Now we use the assumption that $\mathcal{A}\subseteq \{\text{positive integers $\le Q$}\}$ to get a bound on $|\mathcal{B}_{k,s,t}(\mathcal{A},\ell)|$. This gives the inflationary factor. Indeed,
$$ |\mathcal{B}_{k,s,t}(\mathcal{A},\ell)|\ll \prod_{j=1}^{t-1} Q^{j}=Q^{\frac1{2}t(t-1)}.$$
So we have 
$$R_{k,s}(\mathcal{A},\ell)^2\ll Q^{\frac1{2}t(t-1)} \sum_{(\ell_1,\ldots,\ell_{t-1})\in\mathcal{B}_{k,s,t}(\mathcal{A},\ell)}V_{k,s,t}(\mathcal{A},\ell;\ell_1,\ldots,\ell_{t-1})^2.$$
Now we sum over all $\ell$ of the form $\ell=n_1^k+\ldots+n_s^k$ with $n_1,\ldots,n_s\in\mathcal{A}$. This gives
$$\sum_\ell R_{k,s}(\mathcal{A},\ell)^2\ll Q^{\frac1{2}t(t-1)} \cdot \sum_\ell \sum_{(\ell_1,\ldots,\ell_{t-1})\in\mathcal{B}_{k,s,t}(\mathcal{A},\ell)}V_{k,s,t}(\mathcal{A},\ell;\ell_1,\ldots,\ell_{t-1})^2.$$
Finally, we observe what these sums count. The sum on the left-hand side counts the number of solutions to equation (\ref{another_one_equation_system}). The sums on the right-hand side count the number of solutions to our partial Vinogradov system (\ref{partial_vinogradov_system}).

\begingroup
\emergencystretch=3em
\printbibliography
\endgroup

\Addresses

\end{document}